\documentclass[11pt,noinfoline]{imsart}  % <<==  this line will remove the bottom comment on the tex date :) 

\usepackage{amsmath}
\usepackage{amssymb}
\usepackage{amsthm}
\usepackage{amscd,enumerate}
\usepackage{amsfonts}
\usepackage{graphicx}%
\usepackage{fancyhdr, geometry}
\usepackage{url,  color, verbatim}
\usepackage[normalem]{ulem}
\usepackage[colorlinks]{hyperref} 
\usepackage{mathtools}
\usepackage{tikz}
\usetikzlibrary{matrix}
\usepackage{tabularx}
\usepackage{booktabs}

\numberwithin{equation}{section}

\theoremstyle{plain} %\numberwithin{equation}{section}

\theoremstyle{definition}

\begin{document}
\title{Random Monomial Ideals Macaulay2 Package
}
\runtitle{Random monomial ideals}
\begin{center}
\author{%\sout{Richard Osborn[??],} 
Sonja Petrovi\'c, Despina Stasi, Dane Wilburne}
\runauthor{Petrovi\'c, Stasi, Wilburne} 

\address{Department of Applied Mathematics, Illinois Institute of Technology}
%Sonja.Petrovic@iit.edu : Department of Applied Mathematics, Illinois Institute of Technology, Chicago, IL 60616, USA.
%stasdes@iit.edu : Department of Applied Mathematics, Illinois Institute of Technology, Chicago, IL 60616, USA.
%dwilburn@hawk.iit.edu : Department of Applied Mathematics, Illinois Institute of Technology, Chicago, IL 60616, USA.

\end{center}
\maketitle

{\bf Abstract:} 
The {\tt Macaulay2} package {\tt RandomMonomialIdeals} provides users with a set of tools that allow for the systematic generation and study of random monomial ideals. % \cite{RMI}. 
 It also introduces new objects, {\tt Sample} and {\tt Model}, to allow for streamlined handling of random objects and their statistics in {\tt Macaulay2}.  

\medskip
%To do/things to say:
%\begin{itemize}
%\item experiments confirming theorems (PAR 2) <------ added dimStats example
%\item why other RMI does not suffice (tuning parameters, principled study, conjecture generating/testing device, search for counterexamples) (PAR 2)
%\item More paragraph 2: within the framework of a model, you can ask questions about distribution of algebraic invariants of interest <------ Added a sentence about this.
%\item describe what randomMonomialSets does and its overloaded variants (compare/contrast with models in our paper) <------Mostly done. -DW
%\item mention that stats functions will work with more general stuff <-----added this in second section
%\end{itemize}

\section*{Erd\H os-R\'enyi  random monomial ideals}

Given their central role in commutative algebra and their inherent combinatorial structure (see, e.g., \cite{MillerSturmfels, Stanley}), monomial ideals are a natural class of object to study probabilistically. 
This study was initiated in \cite{RMI}, where random monomial ideals were produced from random sets of monomial generators in a manner inspired by the Erd\H os-R\'enyi model of random graphs.  Working within the framework of such a model allows one to ask well-posed questions about the distributions of various algebraic invariants of interest.

The {\tt Macaulay2} \cite{m2} package {\tt RandomMonomialIdeals} implements several basic probabilistic models for monomial ideals based on the work in  \cite{RMI}; it also computes summary statistics of (algebraic properties of) samples of monomial ideals, and sets up the framework %\sout{ to extend the functionality by defining} 
 to define new probabilistic models and generate samples from them using two new {\tt Macaulay2} Types, {\tt Sample} and {\tt Model}. 

\smallskip 
The fundamental method in the % {\tt RandomMonomialIdeals} 
 package is {\tt randomMonomialSets}, which randomly generates sets of monomials in a fixed number of variables up to a given degree from the Erd\H os-R\'enyi-type distribution $\mathcal{B}(n,D,p)$ defined in \cite{RMI}, as well as various other related distributions.  
 % {\color{red} would like to propose editing the revision from here on, to restructure as referee suggested; so what follows is new draft text: } 
 %\hrulefill 
 
In the following examples, we set the number of variables to {\tt n=3} and sample size to {\tt N=5}. 

\medskip
 \begin{tabular}{l|l}
\toprule
 Model & Example of {\tt M2} command \\
\midrule
\midrule
 $\mathcal{B}(n,D,p)$: select each monomial of degree $\leq D$   & {\tt D=3; p=0.2;} \\ 
 %each monomial of degree at most $D$ is selected 
\qquad independently with probability $p\in[0,1]$ &  {\tt randomMonomialSets(n,D,p,N)}\\
\midrule
 $\mathcal{B}(n,D,p)$: select each monomial of degree $1\leq d\leq D$   & {\tt D=3;} \\ 
 \qquad independently with probability $p_d\in[0,1]$, where &  {\tt p = \{0.5,0.0,0.1\};} \\ 
 \qquad $p=(p_1,\dots,p_D)$ is a list of probabilities &  {\tt randomMonomialSets(n,D,p,N)}\\
\midrule
 $\mathcal{B}(n,D,M)$: select $M$ monomials of degree $\leq D$   & {\tt D=3; M=2;} \\ 
 \qquad  uniformly at random &  {\tt  randomMonomialSets(n,D,M,N)} \\  
\midrule
 $\mathcal{B}(n,D,M)$: select $M_d$ monomials of degree $d$   & {\tt D=4; M=\{1,0,3,0\};} \\ 
 \qquad  uniformly at random, where $1\leq d\leq D$ &  {\tt  randomMonomialSets(n,D,M,N)} \\  
\bottomrule

 \end{tabular}
 \medskip
 
% \marginpar{\color{red}$B(n,D,M)$ description in table is correct, isn't it?}
 Each command  in the table above returns a {\tt List} of  $5$   sets of randomly generated monomials.  
In case of $\mathcal B(n,D,M)$, if $M$ is  larger than the total number of monomials in $n$ variables of degree at most $D$,  all such monomials are returned.

%%%%%%%%%%%%%%%%%%%%%%%%%
%%%%%%%%%%%%%%%%%%%%%%%%%

%\smallskip
%Finally, the user may choose for the monomial sets generated to be minimal generating sets. The following example demonstrates the {\tt Minimal} strategy for the graded fixed number of monomials model. It can also be employed for {\tt randomMonomialSets(n,D,p,N)} with $p\in[0,1]$ or $\mathbf{p}\in[0,1]^D$.

One can also force the monomial sets to be minimal generating sets as follows. 
\begin{verbatim}
i1 : n=3; D=4; N=4; p={0.5,0.0,1.0,0.0};  
     netList pack(4,randomMonomialSets(n,D,p,N, Strategy=>"Minimal"))
      +------------+------------------------+------------------------+------------+
      |          3 |      3   2       2   3 |      3   2       2   3 |            |
o1 = |{x , x , x }|{x , x , x x , x x , x }|{x , x , x x , x x , x }|{x , x , x }|
      |  1   2   3 |  3   1   1 2   1 2   2 |  1   2   2 3   2 3   3 |  1   2   3 |
      +------------+------------------------+------------------------+------------+
\end{verbatim}
In the above sample of $4$ sets of monomials, there are no monomials of degrees 2 or 4. Each of the variables was selected with probability $0.5$. All of degree-3 monomials were selected, but of course some are not included as they are not minimal generators of the corresponding monomial ideal. 

\smallskip
%\section*{Generating random monomial ideals}
The distributions that we have described above  and, in fact,  any probability distribution on sets of monomials naturally induce  distributions on monomial ideals. The ideal distribution induced by $\mathcal{B}(n,D,p)$ is denoted by $\mathcal{I}(n,D,p)$.    %; for $\mathcal{B}(n,D,p)$ the corresponding distribution on ideals is denoted $\mathcal{I}(n,D,p)$. %We may sample from  the latter directly as follows:
%To generate 
Samples of monomial \emph{ideals} for each of the models in the table above can be generated  as follows: %from these distributions, replace {\tt randomMonomialSets} with {\tt randomMonomialIdeals} in each of the options of the table above. For example: 
\begin{verbatim}
i2 : n=2;D=5;p=0.2;N=3; randomMonomialIdeals(n,D,p,N)
                                         3
o2 = {monomialIdeal(x ), monomialIdeal (x , x x ), monomialIdeal(x x )}
                     2                   1   1 2                  1 2
o2 : List
\end{verbatim}

%Versions of {\tt randomMonomialIdeals} are provided corresponding to each of the distributions on monomial sets we have described as well.

\section*{Summary statistics}

Given a list of {\tt monomialIdeals}, the package contains an array of methods for computing and summarizing various algebraic invariants of the ideals in the list: %Other summary methods compute statistics on 
Krull dimension, degree, projective dimension, Castelnuovo-Mumford regularity, Betti tables and Betti shapes, and proportion of Borel-fixed ideals or Cohen-Macaulay ones. 

 For instance, the method {\tt dimStats} computes the Krull dimension of $k[x_1,\ldots,x_n]/I$ for each monomial ideal $I$ in the list and returns the mean and standard deviation of the sample.  When the optional input ShowTally is set to true, {\tt dimStats} also returns a histogram of the Krull dimensions of the ideals in the list.

\begin{verbatim}
i3 : B = randomMonomialIdeals(3,10,0.01,1000);
i4 : dimStats(B, ShowTally=>true)
o4 = (1.92, .46, Tally{0 => 1  })
                       1 => 146
                       2 => 785
                       3 => 68
o4 : Sequence
\end{verbatim}
In this sample of $N=1000$ monomial ideals from the distribution $\mathcal{I}(3,10,0.01)$, the proportion of ideals with Krull dimension 0,1,2,3 was 0.001, 0.146, 0.785, 0.068, respectively.  By \cite[Theorem 3.2]{RMI}, the probability that a monomial ideal from this distribution has Krull dimension $t$ for $t=0,1,2,3$ is $0.0009, 0.1458, 0.7963, 0.570$, respectively.  By the same theorem, the expected Krull dimension in this case is 1.9094 whereas the observed sample mean in the example is 1.92.

\begin{comment}
{\color{red}do we ened to keep the example with CMstats? $\rightarrow$ } Another such summary method is {\tt CMStats}, which returns the proportion of monomial ideals in a list whose quotient ring is arithmetically Cohen-Macaulay.

\begin{verbatim}
i10 : CMStats(B, Verbose=>true)
351 out of 1000 ideals in the given sample are Cohen-Macaulay.
There are 1000 ideals in this sample. Of those, 68 are the zero ideal.
They are included in the reported count of Cohen-Macaulay quotient rings.
       351
o10 = ----
      1000
o10 : QQ
\end{verbatim}
\end{comment}

Each method  that computes sample statistics of a particular invariant or property of a list monomial ideals  
% for a list of MonomialIdeals 
%in the {\tt RandomMonomialIdeals} package for computing sample statistics of a particular invariant or property of a list monomial ideals 
 can also be applied more generally to any list of algebraic objects for which that invariant or property is defined, whether or not the objects were generated using the ER-model. 
% {\color{red}make sure we make this last point in the documentation; i think we do?}
% {\color{red} or, in fact, as the referee notes, ``This statement is true, but  one can even more generally run a statistics on any sample with any  numerical function." -- but the referee forgets the CMStats example which is why we included it: this is the true/false and we can still get the summary statistics. However, currently, {\tt statistics} wants to run only on numerical input. The referee suggests generalizing to anything that can be run through {\tt Tally}. Is this difficult?  ... trying to work on it.  - yes it is. The reason it fails is not necessarily the betti command but hte fact that statistics works on a Sample object, not a list. So we do first need to generalize the model construct, which we didnt provide before. I am having huge trouble with it now.}

\section*{New types for creating and storing probabilistic models and samples in Macaulay2}

The package comes equipped with the predefined Erd\H os-R\'enyi-type model. For example, let us consider the graded ER-type model with $p=(0.1,0,0.2)$ in $4$ variables and degree bound $D=3$. We store this in an object of class {\tt Model}:  
\begin{verbatim}
i5 : myModel = ER(ZZ/101[a..d],3,{0.1,0.0,0.2})
o5 = Model{Generate => {*Function[RandomMonomialIdeals.m2:168:22-168:46]*}}
            Name => Erdos-Renyi
                            ZZ
            Parameters => (---[a, b, c, d], 3, {.1, 0, .2})
                           101
o5 : Model
\end{verbatim}
It is now easy to obtain  a sample of 1000 monomial ideals from this model. Note that model parameters %-- that is, things that can be changed by the user --
% -- in this case, the polynomial ring, degree bound $D$ and probability parameter vector $p$ -- 
 are also stored with the Sample object for easy access. 
\begin{verbatim}
i6 : time mySample = sample(myModel,1000);
     -- used 2.32541 seconds
i7 : mySample.ModelName
o7 = Erdos-Renyi
i8 : mySample.Parameters
        ZZ
o8 = (---[a, b, c, d], 3, {.1, 0, .2})
       101
o8 : Sequence
i9 : mySample.SampleSize
o9 = 1000
\end{verbatim}
The raw data (i.e., actual sets of monomials without the parameter values etc.) from the sample can be loaded as follows: 
\begin{verbatim}
i10 :  time myIdeals = getData mySample;
     -- used 0.00005 seconds
i11 : myIdeals_0
        2      2     2
o11 = {b d, b*d , c*d }
o11 : List     
\end{verbatim}
To obtain statistics on any algebraic property of interest in the given sample using the Sample object, one simply runs the {\tt statistics} command. 
Let us look at the distribution of Krull dimensions for this sample of $1000$ monomial ideals:
\begin{verbatim}
i12 : time statistics(mySample, dim@@ideal)
     -- used 0.325259 seconds
                         2053
o12 =  HashTable{Mean => ----, StdDev => .654363, Histogram => Tally{0 => 5  }}
                         1000                                        1 => 166
                                                                     2 => 608
                                                                     3 => 213
                                                                     4 => 8
o12 :  HashTable
\end{verbatim}
In addition, the command {\tt writeSample(Sample,String)} can be used to write a sample to disk (the String is the filename), a feature that will be useful when large samples are generated and their statistics take a lot of computational time (e.g., Gr\"obner bases or  free resolutions).  This command creates a folder in which the model and data are stored. The sample can then be read via calls to the {\tt sample(String)} method; for more details, the user is referred to package documentation.

The new types, {\tt Model} and {\tt Sample}, along with the {\tt statistics} function, allow one to define a new way to sample random algebraic objects, store the data as a proper statistical sample, and study their algebraic properties under the probabilistic regime. Here is a simple example of a model  that generates $M$ polynomials in $n$ variables of degree $D$ randomly using {\tt Macaulay2}'s built-in {\tt random} function: 
\begin{verbatim}
i13 : f=(D,n,M)->(R=QQ[x_1..x_n];apply(M,i->random(D,R)))
o13 = f
o13 : FunctionClosure
i14 : myModel = model({2,3,4},f,"rand(D,n,M): 
			        M random polynomials in n variables of degree D")
o14 = Model{Generate => {*Function[RandomMonomialIdeals.m2:107:22-107:37]*}     }
            Name => rand(D,n,M): M random polynomials in n variables of degree D
            Parameters => {2, 3, 4}
o14 : Model
i15 : mySample = sample(myModel,10); 
\end{verbatim}
The last line produces a sample of size 10 from myModel. 
\begin{comment} here is the first of the 10 random sets:
\begin{verbatim}
i84 : s=getData mySample;

i85 : s_0

       7 2   6       1 2   7                 2 
o85 = {-x  + -x x  + -x  + -x x  + 3x x  + 2x ,
       4 1   5 1 2   2 2   8 1 3     2 3     3 
      --------------------------------------------
      8 2   6       7 2   1       6       1 2  8 2
      -x  + -x x  + -x  + -x x  + -x x  + -x , -x 
      9 1   5 1 2   8 2   4 1 3   7 2 3   4 3  9 1
      --------------------------------------------
                 7 2   3              8 2   2  
      + 10x x  + -x  + -x x  + x x  + -x , x  +
           1 2   9 2   4 1 3    2 3   3 3   1  
      --------------------------------------------
      2       8 2   1              1 2
      -x x  + -x  + -x x  + x x  + -x }
      3 1 2   5 2   2 1 3    2 3   5 3

o85 : List
\end{verbatim}
\end{comment} 
One can  use the {\tt statistics} function to generate an ideal from each of the 10 sets in the sample, compute the Gr\"obner basis, and report its size: 
\begin{verbatim}
i16 : statistics(mySample, numcols@@gens@@gb@@ideal)
o16 = HashTable{Histogram => Tally{6 => 10}}
                Mean => 6
                StdDev => 0
o16 : HashTable
\end{verbatim}
Any function that can be run through {\tt tally} can serve as input function for the {\tt statistics} method; for more extensive examples, the user is directed to the package documentation. 

\section*{Acknowledgements}
The following  undergraduate students   made valuable contributions to the development of the package during summer of 2017: 
Genevieve Hummel, Parker Joncus, Daniel Kosmas, Richard Osborn, Monica Yun, and Tanner Zielinski. The project is funded by the  NSF grant NSF-DMS-1522662, % Randomized and Structure-Based Algorithms in Commutative Algebra,  
the Illinois Institute of Technology College of Science student summer research stipend, and the McMorris student summer research stipend from the Department of Applied Mathematics at Illinois Institute of Technology.

\bibliographystyle {acm} 
\bibliography{RMIpackage}
\end{document}